\documentclass[12pt]{article}
\usepackage{amsmath}
\usepackage{amssymb}
\usepackage{url}

\usepackage[T1]{fontenc}

\usepackage{newcent}

\DeclareMathOperator \dist {dist}
\DeclareMathOperator \heron {Heron}
\DeclareMathOperator \gram {Gram}

\newtheorem{lemma}{Lemma}[section]
\newtheorem{prop}[lemma]{Proposition}
\newtheorem{thm}[lemma]{Theorem}

\newtheorem{defn}[lemma]{Definition}

\newcommand{\E}{{\mathcal E}}

\newcommand{\sqd}{square-density\ }
\newcommand{\sqds}{square-densities\ }

\newcommand{\pp}{parallelepipedum }

\newcommand{\imes}{\times}
\title{Square-densities, and volume forms}
\author{Anders Kock}
\date{}

\begin{document}

\maketitle

\section*{Introduction}

The Greek geometers (Heron et al.) discovered a remarkable formula, 
expressing the area of a triangle   
in terms of the lengths of the three sides. Here, length and area are 
seen as non-negative numbers, which involves,  
in modern terms, formation of {\em absolute value} and {\em square root}.  To 
express the notions and results involved without these non-smooth 
constructions, one can  express the Heron Theorem in terms of the 
{\em squares} of the quantities in question: if $g(A,B)$ denotes the 
{\em square} of the length of the line segment given by $A$ and $B$, the 
Heron formula says that the square of the area of the triangle $ABC$ 
may be calculated by a simple algebraic formula out the three numbers $g(A,B)$, $g(A,C)$, and 
$g(B,C)$. Explicitly, the formula appears in (\ref{zero}) below. In 
modern terms, the formula is  (except for a 
combinatorial constant $-16^{-1}$) the 
determinant of a certain symmetric $4\times 4$ matrix constructed out 
of three  numbers; see (\ref{CMx}) below. This 
determinant, called the Cayley-Menger determinant, 
generalizes to simplices of higher dimesions, so that e.g.\ the square of the 
volume of a tetrahedron (3-simplex) $(ABCD)$ in space is given (except for a 
combinatorial constant) by the determinant 
of a certain $5\times 5$ matrix constructed out of the six square 
lengths of the edges of the tetrahedron (by a formula already known 
in the Renaissance). 

The Heron formula has the advantage that it is symmetric w.r.to 
permutations of the $k+1$ vertices of  a $k$-simplex. Also, 
it  does not  refer to the vector space or  affine  structure of the 
ambient space.

We shall in particular consider the case where the  space, in which the  $k$-simplex 
lives, is a Euclidean space: an affine space $E$ whose associated vector space $V$ is provided 
with a positive definite inner product. Then the square lengths, square areas, square 
volumes etc.\ of the simplices can also be calculated  
by another well known and 
simple expression: namely as $(1/k!)^{2}$ times  the Gram determinant of a certain 
$k\times k$ matrix constructed from the simplex, by choosing one of its vertices 
as origin. The Gram determinant itself 
expresses the square volume of the \pp spanned by $k$ vectors in $V$ 
that go from the origin to the remaining vertices.

An important difference beween the two formulae is the $(k+1)!$-fold 
symmetry in the Heron formula, where the Gram formula is apriori only $k!$-fold 
symmetric, because of the  special role of the chosen origin.

This Gram method of calculating the square-volumes has the advantage that it 
it is easy to describe algebraically, and in particular, it is easy to 
describe what happens if one changes the 
metric; this is needed, when dealing with Riemannian manifolds, where 
 the metric tensor, in any given coordinate chart, changes from point to point.

\medskip

We begin in Section \ref{Sec1x} by recalling the 
classical case of Euclidean spaces. In particular, we recall the comparison 
(standard, but non-trivial) between 
the Heron and Gram calculations. This Section is essentially a piece 
of standard linear algebra. 

In Section \ref{Sec2x}, we recall or 
introduce the notions of differential form and square density in the 
combinatorial versions from synthetic differential geometry (SDG). 
This leads to synthetic, or combinatorial 
arguments, based  on ``infinitesimal'' simplices, and their 
square volume.

In Section \ref{Sec3x}, we relate (in terms of SDG) the volume form of 
an $n$-dimensional Riemannian 
manifold to the volume of certain infinitesimal $n$-simplices. 
 
This Section contains the main theorem, where we, for a 
Riemannian manifold of dimension $n$, compare the square volume of $n$-simplices 
given, respectively, by the Heron 
formula and by the (valuewise) square of the volume form.

 Throughout, $R$ denotes ``the'' number line, a commutative ring with 
 suitable properties to be described when needed. In particular, the 
 notion of {\em positivity}, and of when a quadratic form over $R$ is {\em positive definite} 
 is  recalled in 
 the beginning of Section \ref{S4x}. I do not know whether positive 
 definiteness plays a role in the algebraic arguments in the first 
 three Sections, except that the use of the phrases ``square length'', 
 \ldots , ``square volume'', etc.\ are somewhat misleading in the 
 indefinite case.

It is useful to think in terms of the quantities occurring  as being  
quantities whose physical dimension is some power of length 
(measured in meter $m$, say), so that length is measured in $m$, area 
in $m^{2}$, square area in $m^{4}$, etc.  Tangent vectors are not 
used in the following; they would have physical dimension of $m\cdot t^{-1}$ 
(velocity). The 
word {\em \sqd } is used in any dimension. Square length, square 
area, and square volume are examples, but we do not claim  that the 
square densities considered presently  have such geometric significance.

The theory developed here 
was also attempted in my \cite{[VF]}; I hope that the present account 
will be less ad hoc.

\section{Square volumes in Euclidean spaces}\label{Sec1x}

\subsection{Heron's formula}

The basic idea for the construction of a square $k$-volume function
goes, for the case $k=2$, back to Heron of Alexandria (perhaps even
to Archi\-medes); they knew how to express the square of the area of a 
triangle $S$ (whether located in Euclidean 2-space  or in a higher 
dimensional Euclidean space) in terms of an expression involving only the lengths 
$a,b,c$ of the three sides:
$$area^2 (S)  = t\cdot (t-a)\cdot (t-b) \cdot (t-c)$$
where $t= \tfrac{1}{2} (a+b+c)$. Substituting for $t$, and multiplying out,
one discovers 
(cf.\  \cite{Cox} 1.53) 
that all terms involving an odd
number of any of the variables
$a,b,c$ cancel, and we are left with
\begin{equation}area ^2 (S) = -16^{-1}(a^4 +b^4 +c^4 - 2a^2 b^2 
-  2 a^2
c^2 - 2 b^2
c^2 ),\label{zero} \end{equation}
 an expression that only involves the {\em
squares} $a^2$, $b^2$ and $c^2$ of the lengths of the sides.

The expression in the parenthesis here may be written in terms of the 
determinant of  $4\times 4$ 
matrix (described in (\ref{CMx}) below), which makes it is possible to 
generalize from 2-simplices (= triangles) to $k$-simplices,  in terms of determinants  of certain $(k+2)\times 
(k+2)$ matrices,  ``Cayley-Menger matrices/determinants''; they again only  
involve the square lengths 
of the $\binom{k+1}{2}$ edges of the simplex.

An  $k$-simplex $X$ in a space $M$ is a $(k+1)$-tuple of points 
(vertices) $(x_{0}, x_{1}, \ldots ,x_{k})$ in $M$. If $g: M\times M \to R$ 
satisfies $g(x,x)=0$ and $g(x,y)=g(y,x)$ for all  $x$ and $y$ (like 
a metric $\dist (x,y)$, or its square), one may construct a $(k+2)\imes (k+2)$ 
matrix $C(X)$ by 
the following recipe: first take the  $(k+1)\times (k+1)$ matrix 
whose $ij$th entry is $g(x_{i},x_{j})$. It  has $0$s 
down the diagonal and is symmetric, by the two assumption about $g$. Enlarge this 
matrix  it to a $(k+2)\times (k+2)$- matrix by bordering 
it with  $(0,1, \ldots ,1)$ on the top and  and on the left. The case 
$k=2$ is depicted here (writing $g(ij)$ for $g(x_{i},x_{j})$ for 
brevity; note $g(01)=g(10)$ etc., so that the matrix is symmetric. 
\begin{equation}\label{CMx}\left[ \begin{array}{cccc c}
0&1&1&1\\
1&0&g(01)&g(02)\\
1&g(10)&0&g(12)\\
1&g(20)&g(21)&0\end{array}\right]\end{equation}
(The indices 
of the rows and columns are most conveniently taken to be $-1,0,1,2$.)

This is the Cayley-Menger matrix $C$ for the simplex, and its 
determinant is its 
Cayley-Menger determinant. Heron's formula then says that the value of this determinant is, 
modulo the ``combinatorial'' factor $-16^{-1}$, the square of the {\em area} of a triangle with 
vertices $x_{0}, x_{1}, x_{2}$, as expressed in terms of squares $g(x_{i},x_{j})$ of the {\em 
distances} between them. Similarly for 
(square-) volumes of higher dimensional simplices. Note that no 
coordinates are used in the construction of this matrix/determinant. 

The general formula is that the 
square of the volume of a $k$-simplex is $-(-2)^{-k}\cdot (k!)^{-2}$ times 
the determinant of $C$, e.g.\ for $k=1, 2$, and $3$,  the factors are 
$2^{-1}$,  $-16^{-1}$, and 
 $288^{-1}$, respectively. 
 
 We shall in the following denote the square volume of a $k$-simplex 
$X$, as calculated by the Heron-Cayley-Menger formula, by 
$\heron (X)$ (provided of course that we have some data giving us the 
``square distance'' $g(x_{i},x_{j})$ between its vertices).

 \begin{prop}The Cayley-Menger determinant for a $k$-simplex is invariant 
 under the $(k+1)!$ symmetries of the vertices of the simplex.
 \end{prop}
 {\bf Proof.} Interchanging the vertices $x_{i}$ and $x_{j}$ has the 
 effect of first interchanging the $i$th and $j$th column, and then 
 interchanging the $i$th and $j$th row of the new matrix. Each of these 
 changes will change the determinant by a factor $-1$.
   
\subsection{Gram's formula}

Given a $k$-simplex $X=(x_{0}, x_{1}, \ldots ,x_{k})$ 
in a Euclidean space $E$ with associated vector space $V$ with an 
inner product. If $V=R^{n}$ with the standard inner product, we may 
form an $n\times k$  matrix $Y$
with columns $y_{i}:=x_{i}-x_{0}$ ($i=1, \ldots ,k$). The Gram matrix of the simplex $X$ is then 
the $k\times k$ matrix $Y^{T}\cdot Y$. Its determinant is the Gram 
determinant of the simplex: $\gram(X):= \det( Y^{T}\cdot Y )$. The 
determinant itself is coordiante independent, i.e.\ it only depends 
on the inner product on $V$, not on coordinatizing $V$ by $R^{n}$.

This determinant  likewise has a volume theoretic 
significance: it gives the square of the volume of the 
parallelepipedum spanned by the $k$ vectors $y_{i}:=x_{i}-x_{0}$ in 
$V$ ($i=1, 
\ldots ,k$).

 The following Proposition is only included for a comparison with the 
 issue of $(k+1)!$ symmetry of the formulae.
 
  \begin{prop}\label{GrSx}The Gram determinant for a $k$-simplex is invariant under the $(k+1)!$ 
 symmetries of the $k$-simplex. 
 \end{prop}
 {\bf Proof.} It suffices to prove this for the case where $V=R^{n}$ 
 with standard inner product. Interchanging $x_{i}$ and $x_{j}$, for $i$ and $j\geq 1$ 
 implies an interchange of the corresponding columns in the 
 $Y$-matrix, and this interchanged matrix comes about by multiplying 
 $Y$ on the right by the $k\times k$ matrix $S$   obtained from the unit matrix by 
 interchanging its $i$th and $j$th column. This $S$ has determinant $-1$. 
 So $S^{T}\cdot Y^{T}\cdot Y\cdot  S$ has the same determinant as 
 \mbox{$Y^{T}\cdot Y$}. Interchanging 
 $x_{0}$ and $x_{j}$ in the simplex 
 corresponds, using  $y_{i}=x_{i}-x_{0}$,  to multiplying the $Y$-matrix on the right by the matrix 
 $S_{j}$, 
 obtained from the unit $k\times k$ matrix by replacing its $j$th row 
 by the row $(-1, -1, \ldots,-1)$. The matrix $S$ has determinant $-1$, 
 so $S^{T}\cdot Y^{T}\cdot Y\cdot  S$ has same determinant as 
 $Y^{T}\cdot Y$. Similarly for 
 exchanging $x_{0}$  by $x_{j}$.

\subsection{Comparison formula}
For a Euclidean space $E$, it 
makes sense to compare the values of the Heron and Gram formulas for square volume 
of a $k$-simplex $X=(x_{0}, x_{1}, \ldots ,x_{k})$. 
Let $C$ denote the  $(k+2)\times (k+2)$ matrix ((Heron-) Cayley-Menger)  formed by 
 the square distances between the vertices, as 
described above, and let $Y^{T}\cdot Y$ be the Gram $k\times k$ 
matrix of the simplex, likewise described above. There is a known relation between 
their determinants
 \begin{equation}\label{stackx} \det (C) = - (-2)^{k}\det (Y^{T}\cdot 
 Y).
 \end{equation}
 For a proof, see reference \cite{[Stack]}.

 Note that 
 the left hand hand side in (\ref{stackx}) does not make use of the algebraic structure 
 of $E$ and its associated vector space, but only on the (square-) 
 distance function (arising from the inner product). This flexibility 
 will be crucial when we consider Riemannian manifolds.
 
 We denote the square volume of a simplex $X$, as calculated in terms of 
 the Cayley-Menger matrix $C$, by $\heron (X)$, and  denote the 
 square volume of the corresponding parallelepipedum, as calculated 
 by Gram's method, by 
 $\gram (X)$. But we shall later have occasion to consider different 
 fixed 
 (positive definite) inner products $G$ on 
 one and the same vector space $V$, in which case we may extend the 
 notation and write $\heron_{G}$ and $\gram_{G}$ to specify which 
 inner product we use. The comparison (\ref{stackx}) may then be 
 formulated 
 
 \begin{prop}\label{31x}Let $X=(x_{0}, \ldots ,x_{k})$ be a  $k$-simplex in 
a Euclidean space. Then
$ \heron_{G} (X) = \frac{1}{k!^{2}} 
\gram_{G}(X)$. 
 \end{prop}
 (The factor $k!^{2}$ is just because the volume of the 
 parallelepipedum is $k!$ as large as the one of the simplex itself.)

 \medskip
 
 \noindent{\bf Remark.} In terms of physical dimension alluded to in the Introduction: 
 volume of a $k$-simplex has 
dimension  $m^{k}$, so its square volume has dimension 
$(m^{k})^{2}$; the entries $g(x_{i},x_{j})$ in the Cayley-Menger 
matrix have physical dimension $m^{2}$, and expanding its 
determinant, all terms are products of $k$ copies of these entries. 
(The entries $0$ and $1$ in the top line and left column in the 
matrix are ``pure'' quantities, i.e.\ of dimension $m^{0}$). So the value 
of the determinant is of physical dimension $(m^{2})^{k}$. The Heron 
formula is then meaningful in the sense that it equates quantities of 
 dimension $(m^{2})^{k}$ and $(m^{k})^{2}$.

In particular, the comparison between the square volumes of a 
$k$-simplex, as calculated by Heron-Cayley-Menger and by Gram, which 
is a consequence of (\ref{stackx}), is dimensionally meaningful; both 
have physical dimension $m^{2k}$.

\section{Differential forms and square 
densities}\label{Sec2x} As in \cite{[SGM]}, say, we consider the 
following kind of structure on an object $M$ in a category $\E$ with 
finite inverse limits $M$,  namely 
subobjects of $M\times M$,
$$M_{(0)}\subseteq M_{(1)}\subseteq M_{(2)}\ldots \subseteq M\imes M,
$$
each of the $M_{(r)}$s being a reflexive and symmetric relation, with 
$M_{(0)}$ being the equality relation. We have in mind the ``$r$th 
neighbourhood of the diagonal'' of  an affine 
scheme, as considered in algebraic 
geometry,   or the ``prolongation spaces'' of manifolds as considered 
in e.g.\ \cite{KS}. Except for $M_{(0)}$, these relations are not 
transitive. - We are actually only interested in the the cases 
$r=0,1,2$.

We use the well known ``synthetic'' language to express constructions 
in categories $\E$ with finite limits, in ``elementwise'' 
terms\footnote{Recall that a generalized element of an object $M$ in 
a category $\E$ is just an arbitrary map in $\E$ with codomain $M$; 
see \cite{[SDG]} II.1, \cite{MM} V.5, or \cite{McL} 1.4.}, in 
particular we consider, for a 
natural number $k$, the object of $r$-infinite\-si\-mal $k$-simplices in $M$, meaning the 
subobject of $M \times M \times \ldots \times M$ ($k+1 $ times) 
consisting of $k+1$-tuples $(x_{0}, x_{1}, \ldots  , x_{k})$ of 
elements of $M$ with 
$(x_{i},x_{j}) \in M_{(r)}$ for  all $i,j = 0, 1, \ldots ,k$; such a 
$k+1$-tuple, we shall call an {\em $r$-infinitesimal $k$-simplex}; 
the $x_{i}$s are the {\em vertices} of the simplex.

Note that the question of whether a $k$-simplex is $r$-infinitesimal 
only depends on the ``edges'' $(x_{i},x_{j})$ (face-$1$-simplices) of 
the simplex, equivalently, it depends on the 1-skeleton of the 
simplex.

We shall, as in \cite{[SGM]}, write $x_{i}\sim _{r} x_{j}$ for $(x_{i},x_{j})\in 
M_{(r)}$.  
In the context of SDG, we have that  $x\sim_{r}y$ in $R^{n}$ 
is equivalent to:

{\em For any  $r+1$-linear function $\phi: 
R^{n}\times \ldots \times R^{n}\to R$, we have }\begin{equation}\label{phix}\phi(x-y, \ldots 
,x-y)=0. \end{equation}

\medskip

For $r=1$ and $r=2$, we shall consider certain maps from the object of 
$r$-infinitesimal $k$-simplices to $R$, namely maps which have the property that 
they vanish if $x_{i}=x_{j}$ for some $i 
\neq j$. For $r=1$, combinatorial differential $k$ forms $\omega$ 
have this property. 
(In the context of SDG, such maps are automatically alternating with 
respect to the $(k+1)!$ permutations of the $x_{i}$s, see \cite{[SGM]} 
Theorem 3.1.5.)  

For $r=2$, such maps have not been considered much\footnote{For $r=2$ 
and $k=1$, such things were in \cite{[SGM]} 8.1 called ``quadratic 
differential forms''.}, except for 
the case where $k=1$, where (pseudo-) Riemannian metrics $g$, in the combinatorial 
sense (recalled after Definition \ref{sqdxx} below), are examples of such maps; for this case, we think of 
$g(x_{0},x_{1})$ as the square of the distance between $x_{0}$ and 
$x_{1}$. The $g$s of interest  are  symmetric, 
$g(x_{0},x_{1})=g(x_{1},x_{0})$. For manifolds $M$, we have

\begin{prop}\label{ABx}Given $g:M_{(2)}\to R$ with $g(x,x)=0$ for all $x$. Then 
$g$ is symmetric iff it vanishes on $M_{(1)}\subseteq M_{(2)}$.
\end{prop}
{\bf Proof.} It suffices to consider an $R^{n}$ chart around $x$; we 
consider the 
degree $\leq 2$ part of the Taylor expansion of $g$ around $x$. Then $g$ is
given as $g(x,y)= C(x) +\Omega (x;x-y) + (x-y)^{T}\cdot G(x)\cdot (x-y)$, 
where $C(x)$ is a constant,  $\Omega$ is linear in the argument after the semicolon, 
and $G(x)$ is a symmetric $n\times n$ matrix .
To say that $g$ vanishes on the diagonal $M_{(0)}$ (i.e.\ $g(x,x)=0$ 
for all $x$) is equivalent to 
saying that $C(x)=0$ for all $x$.
We now compare $g(x,y)$ and $g(y,x)$; we claim 
\begin{equation}\label{abx}(x-y)^{T}\cdot G(x)\cdot (x-y) = 
(y-x)^{T}\cdot G(y)\cdot (y-x).\end{equation}
Taylor expanding from $x$ the $G(y)$ on the right hand side gives that this the difference 
between the two sides is $(y-x)\cdot dG(x;y-x)\cdot (y-x)$ 
which is trilinear in $y-x$, and therefore vanishes, since 
$x\sim_{2}y$.
So we have that if $\Omega $ vanishes, then $g$ is symmetric; vice 
versa, if $g$ is symmetric, 
its restriction to $M_{(1)} $ is likewise symmetric, and (being a 
differential 1-form), it is alternating, so the $\Omega$-part 
vanishes, which in coordinate free terms says: $g(x,y)=0$ for $x\sim 
_{1}y$.  

\medskip

For the number line $R$, $(x_{0}, x_{1})\in R_{(2)}$ iff 
$(x_{0}-x_{1})^{3}=0$, and the map $g$ 
given by $g(x_{0},x_{1}):=(x_{0}-x_{1})^{2}$ is  a  map as described in the 
Proposition. In fact, it 
is the restriction of the standard ``square-distance'' function 
$R\imes R \to R$.

So we recall, respectively pose, the following definitions, 
corresponding to $r=1$ and $r=2$. Let $M$ be a manifold.

\begin{defn}A  (combinatorial) {\em differential $k$-form} on $M$ is  an 
$R$-valued   function 
$\omega$ on the set of 1-infinitesimal $k$-simplices in $M$, 
which is alternating
with respect to the $(k+1)!$ permutations of the vertices of the 
simplex.  
\end{defn}
Hence it vanishes on simplices where two vertices are equal.

\begin{defn}\label{sqdxx}A {\em $k$-\sqd } on $M$ is an $R$-valued  function  
on the set of   $2$-infinitesimal $k$-simplices in $M$, 
which is symmetric 
with respect to the $(k+1)!$ permutations of the vertices of the simplex, and which 
vanishes on simplices where two vertices are equal. 
\end{defn}

Note that for $k=1$, Proposition \ref{ABx} gives that
1-square densities (square lengths) $g$ have the property that 
they vanish not just on $M_{(0)}$ (the diagonal), but also on 
$M_{(1)}$: $g(x,y)=0$ if $x\sim _{1}y$. So the notion of 1-square 
density agrees with   (combinatorial) 
``differential quadratic form'', as considered in \cite{[SGM]}, Section 8.1. 
(Combinatorial) differential quadratic 1-forms we shall also call pseudo-Riemannian 
metrics.

\medskip

As a bridge between square densities and differential forms, we pose 
the following auxiliary

\begin{defn}\label{exx}An {\em extended  $k$-form} on $M$ is  an 
$R$-valued   function 
$\overline{\omega}$ on the set of 2-infinitesimal $k$-simplices in $M$, 
which 
vanishes on simplices where two vertices are equal.
\end{defn}

Such extended $k$-form restricts to a function on 1-infinitesimal 
$k$-simplices (and the restriction may or may not be a differential 1-form; note 
that we did not put conditions like ``alternating'' or 
``symmetric'' on extended $k$-forms).

\begin{prop}\label{uniquex}If two extended $k$-forms $\overline{\omega}$ and 
$\overline{\omega}'$ extend the same differential 
$k$-form $\omega$, then  
$\overline{\omega}^{2}=\overline{\omega}'^{2}$.
\end{prop}
{\bf Proof.} We have to prove that  
$$\overline{\omega}^{2}(x_{0}, x_{1}, \ldots ,x_{k})  = 
\overline{\omega}'^{2}(x_{0}, x_{1}, \ldots ,x_{k}),$$
for any 2-infinitesimal 
$k$-simplex $(x_{0}, x_{1}, \ldots ,x_{k})$.
It suffices to do this in a coordinate patch around $x_{0}$, which we 
may assume is $0\in R^{n}$, in which case $\overline{\omega}$ and 
$\overline{\omega}'$ are functions $\Omega$ and $\Omega':  
D_{2}(n)\times \ldots \times D_{2}(n)\to R$ ($k$ factors in the 
product). By the basic axiom scheme of SDG, the ring $A$ of 
functions 
$D_{2}(n) \to R$ is of the form $A=A_{0}\oplus A_{1}\oplus A_{2}$, with 
$A_{0}$ the constant functions $R^{n}\to R$, $A_{1}$ the linear 
functions $R^{n}\to R$, and $A_{2}$ the (homogeneous) quadratic functions $R^{n}\to 
R$. This $A$ is a graded ring (only non-zero in degrees 0,1 and 2). 
The ideal of functions vanishing on $0$ is $A_{1}\oplus 
A_{2}\subseteq A$. So the ideal of functions $(D_{2}(n))^{k} \to R$ 
which vanish if any of its arguments is 0 is the $k$-fold 
(symmetric) tensor product of $(A_{1}\oplus 
A_{2})$,
\begin{equation}\label{idealx}(A_{1}\oplus 
A_{2})^{\otimes k}\subseteq A^{\otimes k}.\end{equation}
This ring is $k$-graded, with e.g.\ the multidegree $(1, \ldots ,1)$ 
consisting of the $k$-linear functions $(R^{n})^{k}\to R$

By assumption, both $\Omega$ and $\Omega'$ belong to the ideal 
(\ref{idealx}).
The assumption that  both $\Omega$ and $\Omega'$ restrict to the same 
differential $k$-form $\omega$ implies that $\Omega$ and $\Omega'$ 
agree in their component of multidegree $(1,\ldots,1)$ (this 
component being the coordinate expression of $\omega$). Thus
$\Omega'=\Omega+\theta$ with $\theta$ of multidegree $\geq (1, \ldots 
,1)$ and of 
total degree $\geq k+1$. The required equation is, in these terms, that 
$(\Omega+\theta)^{2}=\Omega^{2}$, and this is a simple ``counting 
degrees''-argument in the $k$-graded ring $A^{k}$:
\begin{equation}\label{sqx}(\Omega+\theta)^{2}= 
\Omega^{2}+ 2\Omega \cdot \theta + \theta^{2}.\end{equation}
Here, $\theta^{2}$ has total degree $\geq 2\cdot (k+1)\geq 2k+1$, 
which is 0 since $A^{k}$ is 0 in total degrees $>2k$; and $\theta$ is 
a linear combination of terms of multidegree 
of the form $(1, 1, \ldots ,1+p, \ldots 1)$ for $p\geq 1$, so 
$\theta \cdot \omega $ is a linear combination of terms of multidegree
$$(1,1,  \ldots , 1+p, \ldots ,1) + (1,1, \ldots ,1, \ldots ,1) =(2,2, \ldots ,2+p, 
\ldots , 2)$$
which is of total degree $2k+p\geq 2k+1$. So the two last terms in 
(\ref{sqx}) are 0, and this proves the Proposition.

\subsection{$k$-square-densities from $1$-\sqds $g$}\label{HCM} 
We shall argue that for 2-infinitesimal simplices $(x_{0}, \ldots 
,x_{k})$, the 
Cayley-Menger determinants  define square-densities. We already argued 
above that these determinants are symmetric: the value does not change when 
interchanging $x_{i}$ and $x_{j}$. We have to argue for the vanishing 
condition required.
If $x_{i}=x_{j}$, then $g(x_{i}, x_{m})= g(x_{j}, x_{m})$ for all 
$m$, and this implies that the $i$th and $j$th rows in the 
Cayley-Menger matrix are 
equal, which implies that the determinant is $0$.

We denote the $k$-\sqd  corresponding to a $1$-\sqd $g$ by
 $\heron_{g}$ (when $k$ is understood from the context).

\subsection{$k$-square-densities  from differential 
$k$-forms }\label{twox}
Essentially this is the process of {\em squaring} (in $R$) the 
values, so it is tempting to denote the \sqd  which we are aiming for,
by $\omega^{2}$. Precisely:  we get  a well defined $k$-\sqd\ out of a 
differential $k$-form by a two step 
procedure: 1) to {\em extend} the given $k$ form  $\omega$ to a 
suitable function 
$\overline{\omega}$, to allow as inputs not just 1-infinitesimal 
$k$-simplices, but also 2-infinitesimal 
$k$-simplices; and then  2) {\em squaring} $\overline{\omega}$ valuewise.
``Suitable'' means  that $\overline{\omega}$ is an extended form in 
the sense of Defintion \ref{exx}, i.e.\ that it vanishes on simplices 
where two vertices are equal.
 
We shall prove that such an  extension $\overline{\omega}$ is possible; it is not 
unique: it depends on 
choosing a coordinate chart. But we shall prove that uniqueness holds 
after squaring.

The question of existence of such $\overline{\omega}$ is local, so let us assume 
that the manifold $M$ is an open subset of $R^{n}$. Then the $k$-form 
$\omega$  is given by a  function $\Omega: M\times (R^{n})^{k}\to R$, such that
$$\omega (x_{0}, x_{1}, \ldots ,x_{k})= \Omega (x_{0}; x_{0}-x_{1}, 
x_{0}-x_{2}, \ldots , x_{0}-x_{k})$$
where for each $x_{0}\in M$, the function $\Omega(x_{0}; -,\ldots 
,-): (R^{n})^{k}\to R$ is $k$-linear and alternating in the $k$ 
 arguments; these arguments are arbitrary vectors in $R^{n}$, in 
particular, they may be of the form $x_{i}-x_{0}$ for 
$x_{i}\sim_{2}x_{0}$, so the restriction of $\Omega (x_{0}; 
x_{1}-x_{0}, x_{2}-x_{0}, \ldots , x_{k}-x_{0})$ to the the set of 
2-infinitesimal 2-simplexes defines an extension 
$\overline{\omega}$ of $\omega$, so
\begin{equation}\label{deffx}\overline{\omega }(x_{0}, x_{1}, \ldots ,x_{k}):= 
\Omega (x_{0}; x_{0}; x_{0}-x_{1}, 
x_{0}-x_{2}, \ldots , x_{0}-x_{k})
\end{equation}
In this form, the fact that $\overline{\omega}$ is alternating w.r.to the $k!$ 
permutations of the $x_{i}$s ($i=1, \ldots ,k$) can be read of from 
the fact that
$\Omega(x_{0}; \ldots )$ is alternating.
It is also alternating w.r.to permutations involving $x_{0}$, as long 
as the $x_{i}$s are $\sim_{1} x_{0}$; this
can be seen from  seen from an easy Taylor expansion argument, see the 
proof of Theorem 3.1.5 in [SGM]. Now if we use $\Omega$ to 
construct the extension of $\omega$ to 
$\overline{\omega}$, defined on 2-infinitesimal $k$-simplices,  the 
constructed $\overline{\omega}$ will still 
be alternating w.r.to permutations of the $x_{i}$s for $i> 0$, but 
the  Taylor expansion argument mentioned fails for the interchange 
of, say,  
$x_{0}$ and $x_{1}$: we cannot conclude that $\Omega (x_{1}; x_{0}-x_{1}, 
 \ldots )=-\Omega(x_{0}; x_{1}-x_{0}, 
 \ldots )$. This failure get repaired by valuewise squaring:
 \begin{prop}For any 2-infinitesimal $k$-simplex $(x_{0},x_{1}, 
 \ldots ,x_{k})$, we have $$\Omega (x_{1}; x_{0}-x_{1},x_{2}-x_{1}, 
 \ldots )^{2} = \Omega(x_{0}; x_{1}-x_{0}, x_{2}-x_{0}, 
 \ldots )^{2}.$$
 \end{prop}
 {\bf Proof.} We shall only do the case $k=1$. 
  (For  the more general case, the further argument is essentially the same as 
 in the proof of Proposition \ref{GrSx} above.) 
 First, we have  
 by a Taylor expansion from $x_{0}$
 $$\Omega (x_{1}; x_{0}-x_{1})= \Omega (x_{0}; x_{0}-x_{1})+d\Omega 
 (x_{0}; x_{1}-x_{0}; x_{0}-x_{1})$$ $$ + \mbox{ a  term 
 $d^{2}\Omega(x_{0}; \ldots)$, trilinear in } 
 x_{1}-x_{0}.$$
 The trilinear term vanishes, because $x_{1}\sim_{2}x_{0}$. Now we 
 square, and get
 $$\Omega (x_{1}; x_{0}-x_{1})^{2}=\Omega (x_{0}; x_{0}-x_{1})^{2}+ 
 2\cdot \Omega (x_{0}; x_{0}-x_{1})\cdot d\Omega (x_{0}; x_{1}-x_{0}, 
 x_{0}-x_{1})$$ $$+ \mbox{ a term $(d\Omega (x_{0};\ldots))^2$, quadrilinear in } 
 x_{1}-x_{0}.$$
 The quadrilinear term vanishes because $x_{1}\sim_{2}x_{0}$, but also 
 the term $\Omega \cdot d\Omega$ vanishes, because it is trilinear in 
 $x_{1}-x_{0}$. So we get
 $$\Omega (x_{1}; x_{0}-x_{1})^{2}=\Omega (x_{0}; 
 x_{0}-x_{1})^{2}=(-\Omega (x_{0}; x_{1}-x_{0}))^{2} = \Omega (x_{0}; x_{1}-x_{0})^{2},$$
 as desired.
 
 \medskip
 
 We conclude  that a differential 
 $k$-form $\omega$ can be extended to an $\overline{\omega}$ (whose 
 input are 2-infinitesimal $k$-simplices), such that 
 $\overline{\omega}^{2}$ is $(k+1)!$- symmetric. (Also, the extension 
 constructed also clearly has the property that it vanishes if 
 $x_{i}=x_{j}$ for some $i\neq j$.) Hence $\overline{\omega}^{2}$ is 
 a square density. 

 From Proposition \ref{uniquex}, we therefore conclude that
if two extended 
 $k$-forms extend the same differential $k$-form $\omega$, the two 
 resulting square-densities agree.
 
 Because of the Proposition, there is a well-defined ``squaring'' 
 process, leading from differential $k$-forms to $k$-\sqds on a 
 manifold $M$: extend 
 the  form $\omega$, and square the result. It is natural to denote 
 this square density by $\omega^{2}$, with the understanding that it 
 means $\overline{\omega}^{2}$ for any extended form 
 $\overline{\omega}$, extending $\omega$.

 \section{Variable metric tensor}\label{Sec3x}
  
We consider a manifold $M$ which is embedded as an open subset of 
$R^{n}$ (elements of $R^{n}$ we write as $n\imes 1$ matrices). 
A 1-\sqd $g$ on $M$ can in this case be given by a metric tensor, 
i.e.\ by a family 
of symmetric $n\times n$ matrices $G(x)$ (for $x\in M$), such that for 
$x\sim_{2}y$,
\begin{equation}\label{recx}g(x,y) = (x-y)^{T}\cdot G(x)\cdot  
(x-y)\end{equation} (which equals 
  $(y-x)^{T}\cdot 
G(y)\cdot  (y-x)$ by (\ref{abx})).

We shall also use the notation $G(x;;x-y) :=(x-y)^{T}\cdot G(x)\cdot  
(x-y)$. Thus $G(x;;-)$ is quadratic in the argument after the double 
semicolon.

 The letter $G$ is used for the ``metric tensor'', 
 i.e.\ for the family of the 
 matrices $G(x)$. 
 So this $G$ 
 suffices to describe a Heron-Cayley-Menger matrix for any 
 2-infinitesimal $k$-simplex in $M$.   We write  $\heron_{G}(X)$ for 
 the determinant of this matrix. This $\heron_{G}$ defines  in fact a 
 $k$-square density on $M$, for any $k$: {\em metric tensors define 
 square densities}.

We shall prove (Proposition \ref{propsamex})  that for a 2-infinitesimal 
 $k$-simplex, $(x_{0},x_{1},\ldots ,x_{k})$, the $G(x_{i})$s occurring in the 
 Cayley-Menger  determinant for this simplex may all be replaced by $x_{0}$, so 
 that, for a given 2-infini\-te\-simal $k$-simplex,  we can use the comparison 
 with the Gram description, available for constant metric 
 tensors.

The terms in the Cayley-Menger determinant for a $k$-simplex 
$X$ are linear combinations of  $k$-fold products 
$g(x_{i},x_{j})$ with $i\neq j$, in particular the product
\begin{equation}\label{sixx}\pm g(x_{0},x_{1})\cdot g(x_{1},x_{2})\cdot \ldots \cdot 
g(x_{k-1},x_{k})\end{equation}
is a term.
(The other terms in the determinant come about from similar $k$-chains of adjacent 
1-simplices, by permutation of the indices.)  

In terms of variable Riemannian tensors $G(x)$ (with the $G(x)$ symmetric 
$n\times n$ matrices), the product (\ref{sixx}) is (possibly modulo sign) the displayed 
expression in the following Lemma \ref{samex2}). It is useful first 
to intoduce some ad hoc terminology.

 A finite sequence of points $x_{0}, x_{1} \ldots , x_{k}$ in $M$ 
 which are consecutive 2-neigbours i.e.\ $x_{i}\sim_{2}x_{i+1}$ for 
 $i=0, \ldots k-1$, we shall for simplicity call {\em path} of length 
 $k$. If $\tilde{x}$ is a path of length $k$, we get a path of length 
 $k-1$ by omitting the first of the vertex of the path. Let us 
 denote this truncated path by $|\tilde{x}$.
 
 We are interested in such paths in $M \subseteq R^{n}$ when $M$ 
 is eqipped with a Riemannian metric $g$, given by variable symmetric 
 $n\times n$ matrices $G(x)$. So $g(x,y)= (x-y)^{T}\cdot G(x) \cdot 
 (x-y)$. Then for a path $x_{0}, \ldots ,x_{k}$, as above, we write 
 $G(\tilde{x})$ for the product $g(x_{0},x_{1})\cdot \ldots \cdot 
 g(x_{k-1}, x_{k})$, i.e.\ in coordinates 
 \begin{equation}\label{prox}G(\tilde{x}):=G(x_{0};; x_{0}-x_{1})\cdot G(x_{1};; x_{1}-x_{2})\cdot  
  \ldots \cdot G(x_{k-1};; x_{k-1}-x_{k}),\end{equation}
 and we write $\overline{G}(\tilde{x})$ for the similar product, but 
 with all the $x_{i}$s appearing before the double semicolon replaced 
 by the $G(x)$ for $x$ the {\em first} vertex of the path, 
 $$\overline{G}(\tilde{x}):=G(x_{0};; x_{0}-x_{1})\cdot G(x_{0};; x_{1}-x_{2})\cdot  
  \ldots \cdot G(x_{0};; x_{k-1}-x_{k})$$
  
 Thus in
 $\overline{G}(|\tilde{x})$, the constant matrix used is $G(x_{1})$
 because the first vertex of $|\tilde{x}$ is $x_{1}$.
 
 \begin{lemma}\label{samex2}For any path $\tilde{x}$, 
 $G(\tilde{x})=\overline{G}(\tilde{x})$.
 \end{lemma}
 {\bf Proof.} By induction of the length $k$ of the path. The 
 assertion is clearly true for $k=1$. Assume that it holds for $k-1$. 
 Then
 $$G(\tilde{x})=G(x_{0};;x_{0}-x_{1})\cdot G(|\tilde{x})=
 G(x_{0};;x_{0}-x_{1})\cdot \overline{G}(|\tilde{x}),$$
 by the induction assumption, used for the path $|\tilde{x}$. Now 
 by definition of $\overline{G}(|\tilde{x})$, the equation continues
 $$= G(x_{0};;x_{0}-x_{1})\cdot   G(x_{1};; x_{1} 
 -x_{2})\cdot \cdot G(x_{1};;
 x_{i}-x_{i+1})\cdot \ldots G(x_{1};; \cdot \cdot
 x_{k-1}-x_{k}).$$
 Now we Taylor expand, for fixed $i$,  the displayed factor $G(x_{1};; 
 x_{i}-x_{i-1})$, from $x_{0}$ in the direction $x_{1}-x_{0}$: we have
 $$G(x_{1};; x_{i}-x_{i+1})= G(x_{0};; x_{i}-x_{i+1})+dG(x_{0}; 
 x_{1}-x_{0}; x_{i}-x_{i+1})+ Q$$
 $Q$ is quadratic in $x_{1}-x_{0}$. But the 
 linear term $dG(x_{0}; x_{1}-x_{0}; x_{i}-x_{i-1})$, as well as the 
 quadratic term $Q$, get annihilated by being mutiplied with $G(x_{0};; 
 x_{1}-x_{0})$, since this factor is linear in the $dG$-term (and, 
 even  more so, $Q$). So  altogether, we have an  
 expression (at least) trilinear in $x_{0}-x_{1}$, and therefore it 
 vanishes since $x_{0}\sim_{2}x_{1}$. Therefore, in the product 
 (\ref{prox}), each factor 
 $G(x_{i};; \ldots )$ may be replaced by $G(x_{0}; \ldots )$, and then 
 we have $\overline{G}(\tilde{x})$.

\medskip
\noindent{\bf Remark.} The argument simpli\-fies for the case of  ``restricted'' 
2-infini\-tesimal $k$-simplices, as considered by \cite{[Bar]}, since 
there one has that each of the individual $g(x_{i},x_{j})$ in a 
simplex $(x_{0}, \ldots ,x_{k})$ may be calculated by using $G(x_{0})$. 
 In preliminary versions of the present note, I only 
considered the restricted simplices, but the value of the 
Heron-Cayley-Menger formula 
on such simplices is  probably not enough for characterizing the 
volume form, which is our aim.

\medskip

From the Lemma, we conclude, for any variable metric tensor $G$:

\begin{prop}\label{propsamex}Given a 2-infinitesimal $k$-simplex $X=(x_{0}, \ldots 
,x_{k})$. Then $\heron_{G}(X)=\heron_{G(x_{0})}(X)$.
\end{prop}

Combining  with the 
comparison in (\ref{stackx}), we get
\begin{prop}\label{Samx} Given a coordinate patch $M\subseteq R^{n}$ and a 
2-infini\-tesimal $k$-simplex $X=(x_{0},x_{1}, \ldots ,x_{k})$ in 
$M$. Then
$\heron_{G}(X)= (k!)^{-2}\cdot \gram _{G(x_{0})}(X)$. \end{prop}

\section{Volume form}\label{S4x} 
The volume form is a differential  $n$-form that may be defined on an 
$n$-dimensional manifold $M$, which is equipped with a 
{\em Riemannian} (not just pseudo-Riemannian) metric $g$ and which is 
{\em oriented}. 

We need to describe these terms. A Riemannian metric $g$ on $M$ is one 
which in local coordinate charts is given by {\em positive definite} 
matrices $G(x)$, in the following sense:  

We need to assume that $R$ is equipped with a subset $P\subseteq R$ (the 
{\em positive} elements), stable under addition and multiplication, and 
such that elements in $P$ are invertible and have unique positive square roots.
(The unique positive square root of $a\in P$ is denoted $\sqrt{a}$.)
We also assume the dichomotomy: for every invertible $x$, either $x$ 
or $-x$ is positive.

To  say that a metric tensor 
$G$ is positive definite is to say: for each $x\in M$,   $\det (G(x))$ 
is invertible for all $x$, and $G(x)= H(x)^{T}\cdot H(x)$ for some  
 $n\times n$ matrix $H(x)$. Then $\det 
(G(x)) = (\det 
(H(x))^{2}$, so $\det(G)$ is positive and therefore has a square root  
(in fact $\pm \det (H)$ will serve).

An {\em orientation form} for an $n$-dimensional manifold $M$ is a 
differential $n$-form $\delta$, so that any differential $n$-form on 
$M$ can be written\footnote{Recall from the theory of combinatorial 
differential forms (\cite{[SGM]} (3.1.7)) that $f\cdot \delta$ is  the $n$-form given by
$(f\cdot \delta )(x_{0}, \ldots ,x_{n}):= f(x_{0}) \cdot \delta (x_{0}, 
\ldots ,x_{n})$.}
$f\cdot \delta$ for a unique $f:M\to  R$. In the manifold 
$M=R^{n}$, determinant-formation is an orientation form. An 
orientation on $M$ is given by an orientation form, and $\delta_{1}$  
and $\delta_{2}$ define the same orientation if $\delta_{2} = f\cdot 
\delta_{1}$ for an $f:M\to P\subseteq R$. An $n$-form $\omega$ is 
{\em positive} if it is $f\cdot \delta$ for some $f:M\to P$.

Recall  from the last lines of Section \ref{Sec2x} the 
notation $\omega^{2}$ for the square $k$-volume 
constructed out of a differential $k$-form $\omega$:

\begin{thm}Assume that $g$ is a Riemannian metric on an oriented
$n$-dimen\-sional manifold $M$. 
Then there exists   on $M$ a unique positive  differential $n$-form 
$\omega$ such that $\heron_{g}$  and $\omega^{2}$ agree on all 
2-infinitesimal $n$-simplices; 
it deserves the name {\em volume form} for $g$.\end{thm}
{\bf Proof.} Since the data and assertions in the statement do not 
depend on the choice of a (positively oriented) coordinate chart, it suffices to prove the 
assertion in such. So assume that $M \subseteq R^{n}$ 
is an open subset (with orientation inherited from the canonical one 
$\det$ on $R^{n}$), 
% We use this coordinate chart to 
% construct the extension $\overline{\omega}$ as described in Subsection 
% \ref{twox},   and 
and $G$ is given in terms of the  positive 
definite $n\times n$  
matrices $G(x)$ (for $x\in M$). For the existence of a volume form: Consider the extended 
$n$-form $\overline{\omega}$, given by the formula 
$$\overline{\omega}(x_{0},x_{1}, \ldots ,x_{n}):= \frac{ \sqrt{\det 
G(x_{0})}}{n!}\cdot \det (x_{1}-x_{0}, \ldots , x_{n}-x_{0})$$
for any 2-infinitesimal $n$-simplex $X=(x_{0}, \ldots ,x_{n})$. 
Let $Y$ denote the $n\times n$ matrix with $x_{i}-x_{0}$ as its $i$th 
column. 
Then squaring the defining equality for $\overline{\omega}$ gives
\begin{equation}\label{maxx}\overline{\omega}^{2}(X) =\frac{\det 
G(x_{0})}{n!^{2}}\cdot (\det Y)^{2}=\frac{1}{n!^{2}} \det 
(Y^{T}\cdot G(x_{0})\cdot Y)\end{equation}
using the product rule for determinants and $\det (Y^{T})= \det 
(Y)$. By definition of $\gram$, the equation continues 
$$=\frac{1}{n!^{2}}\gram_{G(x_{0})} (X)=\heron_{G(x_{0})}(X) = 
\heron_{G}(X),$$ using  
 the Heron-Gram 
comparison Proposition \ref{31x} and 
Proposition \ref{Samx}.  

This proves the existence of the claimed differential $n$-form.

For the uniqueness,  if we have two positive $n$-forms $\omega_{i}=f_{i}\cdot 
\delta$ (i=1,2) with $f_{i}:M \to P\subseteq R$, we get for all 
2-infinitesimal $n$-simplices $X=(x_{0}, \ldots ,x_{n})$ that
$$ ((f_{1}(x_{0}))^{2}\cdot \delta (X)) \cdot \delta (X) = 
\heron_{G}(X)=((f_{2}(x_{0}))^{2}\cdot \delta (X)) \cdot \delta (X),$$
and cancelling successively the two factors $\delta(X)$ (using the 
uniqueness of the $f$s describing $n$-forms $\omega_{i}$ in terms of $\delta$), we 
ultimately arrive at $f_{1}(x_{0})^{2}=f_{2}(x_{0})^{2}$, and since 
$f_{i}(x)\in P\subseteq R$ for all $x\in M$, we conclude from 
uniqueness of positive square roots that $f_{1}(x_{0})=f_{2}(x_{0})$. 
Since this holds for all 2-infinitesimal $n$-simplices $(x_{0}, 
\ldots ,x_{n})$, we conclude that $f_{1}=f_{2}$, proving the 
uniqueness.

\medskip
\small

\noindent Aarhus University, Dept.\ of Math.\

\noindent November 2021

 \noindent \url{kock@math.au.dk}
 
 \bigskip

\end{document}